\documentstyle[12pt]{article}

\makeatletter

\newfont{\Bbb}{msbm10 scaled 1\@ptsize00}
\newcommand{\CC}{\mbox{\Bbb C}}

\newcommand{\RR}{\mbox{\Bbb R}}

\newif\if@fewtab\@fewtabtrue

\def\draftdate{\number\day.\number\month.\number\year\ \ \ \hourmin }
{\count255=\time\divide\count255 by 60
\xdef\hourmin{\number\count255}
\multiply\count255 by-60\advance\count255 by\time
\xdef\hourmin{\hourmin:\ifnum\count255<10 0\fi\the\count255}}
\def\ps@draft{\let\@mkboth\@gobbletwo
    \def\@oddhead{}
    \def\@oddfoot
       {\hbox to 7 cm{$\scriptstyle\bf Draft\ version:\ \draftdate$
       \hfil}\hskip -7cm\hfil\rm\thepage \hfil}
    \def\@evenhead{}\let\@evenfoot\@oddfoot}

\def\label#1{\ifnum\draftcontrol=1
 \global\def\draftnote{\scriptsize\tt #1}\fi
 \@bsphack\if@filesw {\let\thepage\relax
   \def\protect{\noexpand\noexpand\noexpand}%
\xdef\@gtempa{\write\@auxout{\string
      \newlabel{#1}{{\@currentlabel}{\thepage}}}}}\@gtempa
   \if@nobreak \ifvmode\nobreak\fi\fi\fi
  \@esphack}

\def\@eqnnum{\hbox to 3cm{\phantom{\rm(\theequation)} \draftnote
                         \hfil}\hskip -3cm {\rm(\theequation)}}

\def\eqnarray{\def\draftnote{{}}\global\@fewtabtrue
\stepcounter{equation}\let\@currentlabel=\theequation
\global\@eqnswtrue
\global\@eqcnt\z@\tabskip\@centering\let\\=\@eqncr
$$\halign to \displaywidth\bgroup\@eqnsel\hskip\@centering\@eqcnt\z@
  $\displaystyle\tabskip\z@{##}$&\global\@eqcnt\@ne
  \hskip 1\arraycolsep \hfil${##}$\hfil
  &\global\@eqcnt\tw@ \hskip 1\arraycolsep
$\displaystyle\tabskip\z@{##}$
\hfil  \tabskip\@centering&\global\@eqcnt\thr@@\llap{##}\tabskip\z@
\cr}

\def\endeqnarray{\@@eqncr\egroup
      \global\advance\c@equation\m@ne$$\global\@ignoretrue}

\def\@@eqncr{\let\@tempa\relax
    \ifcase\@eqcnt \def\@tempa{& & &}\or \def\@tempa{& &}
      \or \def\@tempa{&}
      \or\def\@tempa{}
\fi\@tempa
\if@eqnsw
\if@fewtab\@eqnnum\fi
\stepcounter{equation}\fi\global
\@eqnswtrue\global\@eqcnt\z@\global\@fewtabtrue\cr}

\@addtoreset{equation}{section}

\def\cases#1{\left\{\,\vcenter{\normalbaselines\m@th
    \ialign{$\displaystyle{##}\hfil$&\quad##\hfil\crcr#1\crcr}}\right.}

\def\ct#1{\ifnum\draftcontrol=1{\tt [#1]}\else{\cite{#1}}\fi}
\def\ctz#1#2{\ifnum\draftcontrol=1{\tt [#1,#2]}\else{\cite[#1]{#2}}\fi}

\def\draftcite#1{\ifnum\draftcontrol=1#1\else{}\fi}

\def\@lbibitem[#1]#2{\item{}\hskip -3cm \hbox to 2cm
{\hfil$\scriptstyle\draftcite{#2}$}\hskip
1cm[\@biblabel{#1}]\if@filesw
     {\def\protect##1{\string ##1\space}\immediate
      \write\@auxout{\string\bibcite{#2}{#1}}}\fi\ignorespaces}

\def\@bibitem#1{\item\hskip -3cm \hbox to 2cm
{\hfil \scriptsize\tt\draftcite{#1}}\hskip 1cm
\if@filesw \immediate\write\@auxout
       {\string\bibcite{#1}{\the\value{\@listctr}}}\fi\ignorespaces}

\makeatother

\def\lb#1{\label{#1}}
\def\lab#1{\ifnum\draftcontrol=1{{\tt [#1]} \lb{#1}}\else{\lb{#1}}\fi}
\def\Eq#1{(\ref{#1})}
\def\theequation{{\thesection.\arabic{equation}}}
\def\[{\begin{eqnarray}}
\def\nn{\nonumber}
\def\non{\nonumber \\ }
\def\]{\end{eqnarray}}

\def\een{\end{enumerate}}
\def\ben{\begin{enumerate}}

\def\a{\alpha}
\def\b{\beta}

\def\g{\gamma}
\def\d{\delta}
\def\D{\Delta}
\def\e{\epsilon}
\def\th{\theta}

\def\e{\epsilon}
\def\l{\lambda}

\newtheorem{theorem}{Theorem}
\newtheorem{proposition}[theorem]{Proposition}
\newtheorem{lemma}[theorem]{Lemma}

\newtheorem{definition}{Definition}

\def\half{\frac{1}{2}}
\def\del{\partial}

\def\qlie#1{{\cal L}_\hhh(#1)}  % q-Lie algebra
\def\qqlie{{\cal L}_q(g)} 

\def\uqg{U_\hhh(g)}
\def\uqqg{U_q(g)}
\def\ot{\otimes}

\def\lieb#1{[#1]}  % classical Lie bracket
\def\wb#1{\wbxxx#1;}  % white Lie bracket
\def\wbxxx#1,#2;{[#1\circ#2]}
\def\bbxxx#1,#2;{[#1\bullet#2]}
\def\t#1{\tilde{#1}}  % tilde
\def\kill#1{\killx#1;}  % Killing form
\def\killx#1,#2;{B(#1,#2)}
\def\killtx#1,#2;{\t{B}\left(#1,#2\right)}
\def\R{R}  % set of roots
\def\cp#1{\cpppp#1;}
\def\cpppp#1,#2;{#1\cdot#2}
\def\Ch{\CC[[\hhh]]}
\def\form#1{\formxx#1;}  % form on root space
\def\formxx#1,#2;{\langle #1,#2\rangle}
\def\dlie{{\cal L}_\hhh(g)}
\def\hhh{h}

\begin{document}

\def\draft{\pagestyle{draft}\thispagestyle{draft}
  \global\def\draftcontrol{1}}
\global\def\draftcontrol{0}
%%%%%%%%%%%%%%%%%%%% switch on/off draft version %%%%%%%%%%%%%%%%%%%%%%%
%\draft

\newpage
\begin{titlepage}
\begin{flushright}
{KCL-TH-95-04}\\
{q-alg/9506017}
\end{flushright}
\vspace{1.5cm}
\begin{center}
{\bf {\Large ON QUANTUM LIE ALGEBRAS }\\
\vspace{4mm}
{\Large AND QUANTUM ROOT SYSTEMS}}\\
\vspace{1.2cm}
{\large Gustav W. Delius}
\footnote{On leave from Department of Physics, Bielefeld University,
Germany} and {\large Andreas H\"uffmann}\\
\vspace{3mm}
Department of Mathematics,
King's College London\\
Strand, London WC2R 2LS, Great Britain\\
{\small e-mail: delius@mth.kcl.ac.uk and aha@mth.kcl.ac.uk}\\
{\small www: http://www.mth.kcl.ac.uk/\~delius}
\vspace{1.6cm}

{ABSTRACT}
\end{center}
\begin{quote}
%Motivated by properties of quantum affine Toda theories we develop an
%analysis of quantum Lie algebras in terms of quantum roots and quantum
%structure constants.

As a natural generalization of ordinary Lie algebras we introduce the concept 
of quantum Lie algebras ${\cal L}_q(g)$.
We define these in terms of certain adjoint submodules of quantized
enveloping algebras $U_q(g)$ endowed with a quantum Lie bracket given
by the quantum adjoint action. The structure constants of these algebras
depend on the quantum deformation parameter $q$ and they go over
into the usual Lie algebras when $q=1$.

The notions of q-conjugation and q-linearity are introduced. q-linear analogues
of the classical antipode and Cartan involution are defined and a  generalised
Killing form, q-linear in the first entry and linear in the second,
is obtained. These structures allow the derivation of symmetries between 
the structure constants of quantum Lie algebras.

The explicitly worked out  examples of $g=sl_3$ and $so_5$
illustrate the results.
\end{quote}
\vfill
\end{titlepage}

\section{Introduction\lab{sectintro}}

Lie algebras and their associated root systems play a pervasive
role in the theory of classical integrable models.
The great breakthrough in the quantization of these models has
been the realization of the importance of the quantized enveloping
algebras associated to these Lie algebras \ct{Jim85,Dri85,Dri86,FRT}.
With the help of these quantized enveloping algebras it has been possible to
derive many exact results for the full quantum theories.

In this paper we will deal not with the quantization of the enveloping
algebras of Lie algebras but with the quantization of the Lie algebras
themselves. Given the fact that most of the properties of classical
integrable models are described by the structure
of Lie algebras rather than their enveloping algebras, it is worthwhile to
attempt to describe the quantum integrable models with quantum Lie
algebras instead of quantized enveloping algebras. In section
\ref{sectphysical} we will describe the particular examples of quantum
integrable theories which motivated this work.

A Lie algebra $g$ is naturally embedded into its universal enveloping
algebra $U(g)$ as a submodule with respect to the adjoint action. The
Lie bracket on $g$ is the restriction of the adjoint action of $U(g)$
to this submodule.

In the quantum case we are given the quantized enveloping algebra
$\uqqg$ and its quantum adjoint action on itself. We study
those submodules of $\uqqg$ which under the quantum adjoint action
transform as the adjoint representation, following a remark in
\ct{Kir90}. We endow these modules with the quantum Lie bracket induced
by the quantum adjoint action. The resulting algebras are not all isomorphic.
But among them there are always distinguished ones which share further
important properties with their classical counterparts and it is these
which we study in detail in this paper. The precise definition of these
quantum Lie algebras is contained in Definition \ref{defqlie}.

There is a different approach to the quantization of Lie algebras present in
the literature. It is based
on the notion of bicovariant differential calculus
on quantum groups \ct{Wor89,Ber90,Ber91,Jur91,Asc93,Sch93a,Sch93b,Sch95}.
The resulting structures are braided Lie algebras as discussed by
Majid \ct{Maj93}. Their shortcoming is that they do not have the same 
dimension as the corresponding
classical Lie algebras except in the case of $g=gl_n$.
For a discussion of this problem see \cite{Sud1}. For the case of
$g=sl_n$ this problem has recently been solved by
Sudbery and Lyubashenko \cite{Sud2}.

This paper is structured as follows. In section 2 we briefly mention the
features of affine Toda quantum field theories which motivated our
search for a quantum deformation of Lie algebras and root systems.
This section is included purely as a motivation. 
Section 3 contains some necessary preliminary material on Lie algebras and
on quantum enveloping algebras. In order to introduce the concept of
quantum Lie algebras we give in section \ref{sectexample}
the very simple example of $\qlie{sl_2}$.
In section \ref{sectstruc}
we give the beginnings of a general study of the structure of quantum
Lie algebras. The standard tools provided by the general structure
of quantum groups are complemented with the notion of q-conjugation.
It is this construction that allows us to exploit a generalisation
of the classical Killing form, defined in section \ref{sectkilling},
to obtain the analogue of the Weyl canonical form of a Lie algebra in
section \ref{sectcanon}. Relations and symmetries of the structure
constants of the quantum Lie algebras in this basis are derived in
section \ref{sectrelations} and the quantum root space is investigated
in section \ref{sectroots}.

Finally the structure constants for the quantum Lie algebras associated
with the Lie algebras $a_2$ $(=sl_2)$ and $c_2$ $(=sp(4)=so(5))$ are
given in section \ref{sectex}. The calculations were done on a computer
using Mathematica \ct{Wol91}. The results were obtained without using
the general results of section \ref{sectstruc} on the structure of 
quantum Lie algebras but are of course found to be in agreement with them. 
By the same methods we have also obtained the explicit results for the 
quantizations of the Lie algebras $a_3=sl_4$ and $g_2$. All the explicit
calculations and results are available in the form of Mathematica 
notebooks at http://www.mth.kcl.ac.uk/\~delius/q-lie.html on the 
World Wide Web.

The straightforward determination of the explicit $q$-dependent structure
constants of quantum Lie algebras $\qqlie$ is extremely tedious. We have 
therefore recently described a general method for obtaining them from the
R-matrix of $U_q(g)$ \ct{Del95}. This method had independently and
earlier been derived in the formalism of differential calculus on
quantum groups, see e.g. \ct{Sch93a}.  However in \ct{Del95} it is 
applied to $g=gl_n$ and $g=sl_n$ for all $n$.
The paper \ct{Del95b} establishes
the existence and uniqueness of the quantum Lie algebras discussed here.

\section{Physical Motivation\lab{sectphysical}}

We want to start by giving the physical motivation which has led us to
undertake the present study of quantum Lie algebras and quantum root
systems. This section is meant
purely as a motivation and is in no way needed in the rest of the paper.

This work has grown out of our desire to understand the exact results
which have been obtained in quantum affine Toda theories. In these theories
it has been possible to obtain the full
quantum mass ratios and the exact S-matrices for the fundamental particles
\ct{Bra90,Del92}.
Furthermore, Dorey \ct{Dor91} has found an elegant description of
these results in terms of properties of the root systems of the
underlying Lie algebras. While this description is exact for the cases
where the affine root system is self-dual, the true quantum results in the
case of
non-self-dual root systems require certain deformations, with the
deformation parameter depending on the product of Planck's constant
and the coupling constant \ct{Del92}.

It is tempting to conjecture that
the systematics of these deformations might be understandable
in terms of the quantum root systems of quantum Lie algebras.
However, a concept of quantum root systems associated to quantum Lie
algebras has, to our knowledge,
never been studied in the literature.

Affine Toda theories are massive integrable two dimensional relativistic
field theories described by the Lagrangian density
\[
L[\phi]=b\left(\del_\mu\phi,\del^\mu\phi\right)+
\frac{m^2}{\b^2} b\left(e^{\b\rm{ad}\,\phi}z_1,z_{-1}\right),
\]
where the bosonic field $\phi$ takes its values in the Cartan subalgebra
of a simple Lie algebra $g$, $m^2$ is a mass scale, $\b$ is the coupling
constant, $b(\,.\,,\,.\,)$ is the Killing form on $g$. The $z_{\pm 1}$
are cyclic elements of $g$ which in a standard notation can be
expressed as
\[
z_1=\sum_{\a\in\bar{\D}}\sqrt{n_\a}\,x_\a,~~~~~~
z_{-1}=\sum_{\a\in\bar{\D}}\sqrt{n_\a}\,x_{-\a}.
\]
where $\bar{\D}$ is the set $\D$ of simple roots extended by the root
$\a_0$ which is minus the highest root
(or the highest short root in the case of twisted Toda theories). The
$n_\a$ are the Kac labels defined so that $n_{\a_0}=1$ and
$\sum_{\a\in\bar{\D}}n_\a \a=0$. The classical
masses of the fields can be read off the Lagrangian and their squares
are found to be
the eigenvalues of the matrix (written in terms of some basis $\{h_i\}$
of the Cartan subalgebra)
\[
M^2_{ij}=\sum_{\a\in\bar{\D}}n_\a\,\a(h_i)\,\a(h_j)
\]
Equivalent characterizations of the squares of the masses is as the
eigenvalues of the adjoint action of $z_1\,z_{-1}$ on the Lie algebra
or as the length squared of the projections of certain roots into the
lowest eigenspace of the Coxeter element of the Weyl group.
Numerically this typically leads to values (slight modifications depend
on the particular Lie algebra $g$)
\[\label{masses}
m_a^2= 8 m^2 \sin^2\frac{a\pi}{2h},~~~~~~~~
h=\sum_{\a\in\bar{\D}}n_\a.
\]
where $a$ is the integer labeling the particle and $h$ is the
(twisted) Coxeter number of $g$.

In the quantum theory these masses receive quantum corrections. However,
when the dust settles, it turns out that the exact quantum masses are still
given by the formula in \Eq{masses} but with the Coxeter number $h$
replaced by a ``quantum'' Coxeter number H. When the set $\bar{\D}$
is self-dual (i.e., if $\forall \a\in\bar{\D}$ also $2\a/\a^2\in\bar{\D}$)
this quantum Coxeter numer is equal to its classical value but in
the non-selfdual case it is coupling constant dependent in
the generic form
\[\label{H}
H=h+c\frac{\b^2\hbar/2\pi}{1+\b^2\hbar/4\pi}
\]
where $c$ depends on the particular Lie algebra.
Will it be possible to find a quantum Lie algebraic explanation for
these mass formulae? In particular, is there a natural definition of a
quantum Coxeter number?

The factorized S-matrices for the fundamental particles of affine
Toda theories have been exactly determined by solving the equations
arising from the bootstrap principle \ct{Bra90,Del92}. Dorey \ct{Dor91}
found that solutions to these very stringent bootstrap equations
could be constructed by using the properties of the root systems
of Lie algebras. These solutions describe the S-matrices of the self-dual
Toda theories. They have the special property that the locations of
the poles do not depend on $\hbar$. In addition to Dorey's solutions
there is another set of solutions in which the pole locations depend
on $\hbar$ through the quantum Coxeter number. These solutions give
the S-matrices of non-self-dual Toda theories. Can the reason for the
existence of these solutions  be understood in terms of the properties
of quantum root systems?

\section{Preliminaries\lab{sectprelim}}

For background on Lie algebras see for example \ct{Sam69}.
Let $g$ be a simple complex Lie algebra of rank $r$,
$R$ the set of non-zero roots and
$\a_1,\a_2,\dots,\a_r$ its simple roots.
Let $b: g\ot g\rightarrow \CC$ be the Killing form.
Choose a basis ${\hat{h}_1,\hat{h}_2,\dots \hat{h}_r}$ for the
Cartan subalgebra ${\cal H}$ so that
$b(\hat{h}_i,h)=\a_i(h)~\forall h\in{\cal H}$.
Choose root vectors $\hat{x}_\a$ so that
$b(\hat{x}_\a,\hat{x}_{-\a})=-1$. Then the Lie bracket relations take
the Weyl canonical form
\[\label{liestruc}
&&\lieb{\hat{h}_i,\hat{x}_\a}=-\lieb{\hat{x}_\a,\hat{h}_i}=
\a(\hat{h}_i)\,\hat{x}_\a,~~~
\lieb{\hat{h}_i,\hat{h}_j}=0\non
&&\lieb{\hat{x}_\a,\hat{x}_{-\a}}=-\hat{h}_\a,~~\mbox{where if }
\a=\sum k_i\a_i \mbox{ then }\hat{h}_\a=\sum k_i \hat{h}_i,\non
&&\lieb{\hat{x}_\a,\hat{x}_\b}=N_{\a,\b}\,\hat{x}_{\a+\b}~~~
\mbox{for }\b\neq -\a \mbox{ and }\a+\b\in R.
\]
The $N_{\a,\b}$ are
real numbers which
%satisfy $N_{-\a,-\b}=N_{\a,\b}$ and which
can be determined entirely in terms of the root system.
The scalar product on the root lattice is defined by
\[
\cp{\a,\b}\equiv b\left(\hat{h}_\a,\hat{h}_\b\right)
=\a(\hat{h}_\b).
\]
The Weyl canonical basis is related to the Chevalley canonical basis by
\[\label{wtc}
x_{\pm\a}=\pm\sqrt{\frac{2}{\cp{\a,\a}}}\,\hat{x}_{\pm\a},~~~~
h_i={\frac{2}{\cp{\a,\a}}}\hat{h}_i
\]
In the Chevalley basis all structure constants are integers.
To generate the Lie algebra it is sufficient to consider the simple
root vectors $x^\pm_i=x_{\pm\a_i}$. The relations are then
\[\label{chevalleyrels}
&&\left[ h_i,h_j \right] = 0,~~~~
\left[ h_i,x_j^\pm \right] = \pm a_{ij} x_j^\pm ,~~~~
\left[x_i^+ ,x_j^- \right]  =  \delta_{ij} h_j,\non
&&\mbox{ad}(x_i^\pm)^{1-a_{ij}}(x_j^\pm)=0~~~\mbox{if }i\neq j.
\]
The last relations are the Serre relations. The adjoint action is
defined by the Lie bracket $ad(x)(y)=[x,y]$ and
$a_{ij}=2\cp{\a_i,\a_j}/\cp{\a_i,\a_i}$ is the
Cartan matrix.

The universal enveloping algebra $U(g)$ is the unital associative algebra over
$\CC$ with generators $x_i^+,\ x_i^-,\ h_i$, $1 \le i \le r$ and relations
\Eq{chevalleyrels} in which the Lie bracket is replaced by the
commutator. The quantized eveloping algebra
$U_\hhh (g)$ is an algebra over $\Ch $, the ring of formal power series
in the indeterminate $\hhh$,
with the same set of generators but with the deformed relations
\footnote{We have found \ct{Cha94} to be a generally reliable reference
on quantum groups. Our $x_i^{\pm}$ are related to the $X_i^\pm$ of
\ct{Cha94} by $x_i^+= k_i^{-1/2} X_i^+$ and $x_i^-= X_i^- k_i^{1/2}$.}
\[\label{uqrel}
&&\left[ h_i,h_j \right] = 0,~~~~~
\left[ h_i,x_j^\pm \right] = \pm a_{ij} x_j^\pm , \non
&&\left[x_i^+ ,x_j^- \right]  =  \delta_{ij}
\ \frac{q_i^{h_i} - q_i^{-h_i}}{q_i- q_i^{-1}},
\]
and the quantum Serre relations
\[\label{serre}
\sum_{k=0}^{1-a_{ij}} (-1)^k
\left[ \begin{array}{c} 1-a_{ij} \\ k \end{array} \right]_{q_i}
(x_i^\pm)^k x_j^\pm (x_i^\pm)^{1-a_{ij}-k} = 0 \qquad i \not= j.
\]
Here $\left[\begin{array}{c}a\\b\end{array}\right]_q$ are the
q-binomial coefficients.
We have defined $q_i = e^{d_i \hhh}$ where $d_i$ are
coprime integers such that $d_i a_{ij}$ is a symmetric matrix.
We will use the notation
 $k_i = q_i^{h_i}$ and then the relations
\Eq{uqrel} take the form
\[\label{krels}
k_i k_j = k_j k_i,\quad k_i x_j^\pm k_i^{-1} = q_i^{\pm a_{ij}} x_j^\pm \quad
\mbox{and} \quad \left[x_i^+ ,x_j^- \right] = \delta_{ij} \
\frac{k_i - k_i^{-1}}{q_i - q_i^{-1}}.
\]
Note the technical point that in this paper we do
not work with some rational form $U_q(g)$ but always with the algebra
$\uqg$ defined over $\Ch $.
Indeed it can be seen from the example of $g=a_2$ that in general
our quantum Lie algebras do not exist in the usual {\it adjoint}
rational form but that one would have to use the {\it simply-connected}
rational form.

The Hopf algebra structure of $\uqg$ is given by the comultiplication
\[
\Delta(h_i) &=& h_i \otimes 1 + 1 \otimes h_i, \\
\Delta(x_i^\pm) &=& x_i^\pm \otimes q_i^{h_i/2} +
 q_i^{-h_i/2} \otimes x_i^\pm,\]
the antipode
\[\label{antipode}
S(h_i)= - h_i, ~~~~
S(x_i^\pm) = - q_i^{\pm 1}\,x_i^\pm ,\]
and the counit
\[\label{counit}
\epsilon(h_i) = \epsilon(x_i^\pm) = 0.
\]
The antipode does not square to the identity but rather
\[\label{defu}
S^2(a)=u\,a\,u^{-1}~~~\mbox{with }u=q^{2h_\rho},
\]
where $q=e^\hhh$ and $h_\rho$ is the element of the Cartan subalgebra
determined by $b(h_\rho,h)=\rho(h)\,\forall h\in{\cal H}$ with $\rho$
being half the sum of the positive roots.

The Cartan involution $\th$ is given by the same formulas as in the
classical case
\[\label{cartinvo}
\th(x_i^\pm)=x_i^\mp,~~~\th(h_i)= -h_i.
\]
It is an algebra automorphism and a coalgebra antiautomorphism
\[\label{coanti}
\D\cdot\th=(\th\ot\th)\cdot\D^T,~~~~~
S\cdot \th=\th\cdot S^{-1}.
\]
If the Dynkin diagram of $g$ has a symmetry $\tau$ which maps node
$i$ into node $\tau(i)$ then the Lie algebra $g$ has an
automorphism
\[
\tau(x^\pm_i)=x^\pm_{\tau(i)},~~~~~\tau(h_i)=h_{\tau(i)}
\]
which extends to a Hopf-algebra automorphism of $\uqg$.
Such $\tau$ are refered to as diagram automorphisms and
except for rescalings of the $x^\pm_i$ they are the only
Hopf-algebra automorphisms of $\uqg$.

The adjoint action of $\uqg$ on itself, using Sweedler's notation \ct{Swe69},
is given by
\[\label{adjoint}
x\circ y=\sum x_{(1)}\,y\,S(x_{(2)}),~~~~~x,y\in\uqg.
\]
There is a second
adjoint action $\bullet$ defined by
\[\label{bad}
x\bullet y=\sum x_{(2)}\,y\,S^{-1}(x_{(1)}).
\]
The Cartan involution $\th$ and the antipode $S$ respect the
adjoint actions in the sense of
$[\th(a)\bullet\th(b)]=\th([a\circ b])$ and
$[S(a)\bullet S(b)]=S([S^{-1}(a)\circ b])$ for all
$a,b\in\uqg$.

\section{The $sl_2$ example}\label{sectexample}

As an introduction to the idea of a quantum Lie algebra it is useful
to consider the very simple example of $sl_2$. The quantized enveloping
algebra $U_\hhh(sl_2)$ is generated by the three generators $h,x^+,x^-$
and the commutation relations
\[
[h,x^\pm]=\pm 2 x^\pm,~~~~~
[x^+,x^-]=\frac{q^h-q^{-h}}{q-q^{-1}}.
\]
Thus these three generators do not close to form a Lie algebra
because the right hand side of the second equation is non-linear.
Of course one would not expect them to do so. In the quantum case
the commutator, which describes the classical adjoint action, should
be replaced by the quantum adjoint action described in \Eq{adjoint}.
In general the adjoint action on any $a\in \uqg$ is given by
\[
h\circ a=[h,a],~~~~~
x^\pm\circ a=x^\pm\,a\,q^{-h/2}-q^{-h/2\pm 1}\,a\,x^\pm.
\]
and this produces the commutator only for $q=1$. The generators
$h,x^\pm$ do not close even under the quantum adjoint action.
However the elements
\[
X^\pm=q^{h/2}x^\pm,~~~~~H=q^{-1}x^+x^--qx^-x^+
\]
do. Indeed, their adjoint actions on each other can be easily
calculated to be given by
\[\label{sl2q}
[H\circ X^+]=(1+q^{-2})X^+,~~~~&&[X^+\circ H]=-(1+q^2)X^+,
\non{}
[H\circ X^-]=-(1+q^{2})X^-,~~~~&&[X^-\circ H]=(1+q^{-2})X^-,
\non{}
[X^+\circ X^-]=H,~~~~~&&[X^-\circ X^+]=-H,
\non{}
[H\circ H]=(q^{-2}-q^2)H,~~~~~&&[X^\pm\circ X^\pm]=0.
\]
We use the bracket notation for the quantum adjoint action to indicate
that we now view it as the quantum analoge of the Lie bracket.
The algebra in \Eq{sl2q} is the quantum Lie algebra $\qlie{sl_2}$. 
Its structure constants
are $q$-dependent in such a way that it goes over into the classical
$sl_2$ Lie algebra for $q=1$.

The simplicity of this example is deceptive. For any Lie algebra
other than $sl_2$ the associated quantum Lie algebra is much more complex.
We give other examples in section \ref{sectex}.

\section{General structure\lab{sectstruc}}

It is now our aim to make some general statements about the structure
of quantum Lie algebras and to derive symmetries between their 
structure constants.

\subsection{q-conjugation \lab{sectqconj}}

%$\uqg$ carries a second coalgebra structure given by the opposite
%coproduct $\D^T$ and the inverse antipode $S^{-1}$. This second coalgebra
%structure can alternatively be viewed as arising from the replacement
%$\hhh\mapsto -\hhh$, i.e. $q\mapsto q^{-1}$.
An important role is played in our general study by the concept of
$q$-conjugation.

\begin{definition}
a) \underline{q-conjugation}
\mbox{$\sim: \Ch \rightarrow\Ch$}, $a\mapsto\t{a}$ is the ring
automorphism defined by $\t{\hhh}=-\hhh$.

\noindent b) Let $M,N$ be $\Ch$-modules. A map
$\phi:M\rightarrow N$ is \underline{$q$-linear} if
\[\label{qlinear}
\phi(\l \,a)=\t{\l}\,\phi(a),~~~~~~~\forall a\in M, \l\in\Ch.
\]
c) Let $A,B$ be algebras over $\Ch$. A q-linear map
$\phi: A\rightarrow B$ is an \underline{algebra $q$-homo-}
\underline{morphism} if it
respects the algebra product, i.e., if
$\forall a,a'\in A,~~\phi(a\,a')=\phi(a)\,\phi(a')$.
$q$-anti-isomorphisms, $q$-automorphisms, etc., are defined analogously.
\end{definition}
Note the analogy between the concepts of $q$-conjugation and complex
conjugation and between $q$-linear maps and anti-linear maps.

\begin{definition}\lab{qu}
\underline{$q$-conjugation} on the quantum group $\uqg$ is the algebra
q-automor\-phism
\mbox{$\sim: \uqg \rightarrow \uqg$} that extends $q$-conjugation on $\Ch$
by acting as the identity on the generators $x_i^\pm$ and $h_i$.
\end{definition}
This definition is consistent because the relations \Eq{uqrel} and
\Eq{serre} are invariant under $q\mapsto q^{-1}$. The notion of
q-conjugation has been introduced already in \ct{Dri90}.

Defining a tilded Cartan involution and a tilded antipode as compositions
\[
\t{S}=\sim\cdot S,~~~~~\t{\th}=\sim\cdot\th,
\]
the concept of q-conjugation proves to be useful as we have

\begin{lemma}
a) $q$-conjugation is a Hopf algebra q-isomorphism
\mbox{$\sim: \uqg \rightarrow \uqg^{\rm{op}}$}, in particular
\[
\e \cdot \sim = \sim \cdot \e,~~~~~
\D \cdot \sim = \sim \cdot \D^T,~~~~~
S \cdot \sim = \sim \cdot S^{-1}.
\]
b) $q$-conjugation relates the adjoint actions as
\[\label{bwr}
\t{a}\bullet\t{b}=\widetilde{a\circ b},~~~~~~~\forall a,b\in\uqg.
\]
c) $\t{\th}:\uqg \to \uqg$ is a Hopf algebra $q$-isomorphism, especially
\[\label{tiso}
\t{\th}(a)\circ\t{\th}(b)=\t{\th}(a\circ b),~~~~~~~\forall a,b\in\uqg.
\]
d) $\t{S}:\uqg \to \uqg$ is an algebra q-anti-isomorphism such that
\[\label{swb}
\t{S}(a)\circ\t{S}(b)=\t{S}(S^{-1}(a)\circ b),~~~~~~~\forall a,b\in\uqg.
\]
\end{lemma}

\subsection{Quantum Lie algebras $\qlie{g}$\lab{sectqlie}}

A Lie algebra $g$ is naturally embedded into its universal enveloping
algebra $U(g)$. It forms a subspace of the enveloping algebra which
under the adjoint action transforms in the adjoint representation and
the adjoint action restricts to the Lie bracket.
As a starting point it is natural to  define a quantum Lie algebra
$\qlie{g}$ as a submodule of the quantized enveloping algebra $\uqg$
with the analogous property. The following definition additionally asks
for a quantum Lie algebra to be invariant under $\t{\th}$, $\t{S}$ 
and $\tau$ as
this is not guaranteed by the classical limit itself.

While a modification of the following definition would be appropriate also
to the case of Kac-Moody algebras, in this paper we have the case
of finite dimensional Lie algebras in mind.

\begin{definition}\lab{defqlie}
A \underline{quantum Lie algebra $\qlie{g}$}
associated to a finite-dimen\-sio\-nal simple
complex Lie algebra $g$ is a finite-dimensional 
indecomposable $\circ$ - submodule
of $\uqg$ endowed with the 
\underline{quantum Lie bracket} 
$[a\circ b]=a\circ b$
such that
\begin{enumerate}
\vspace{-1mm}
\item $\qlie{g}$ is a deformation of $g$, i.e., 
$\qlie{g}|_{\hhh=0}=g$.
\vspace{-1mm}
\item
$\qlie{g}$ is invariant under $\t{\theta}$, $\t{S}$ and any diagram
automorphism $\tau$.
\end{enumerate}
\end{definition}

An immediate consequence of this definition is that,
under the adjoint action of $\uqg$, $\qlie{g}$ transforms
as the adjoint representation. The structure
of this representation is well know. As is the case with all finite
dimensional highest weight representations of $\uqg$ \ct{Ros88,Lus88},
it is just a deformation of the corresponding classical representation. It
follows in particular that $\qlie{g}$ splits into submodules of
definite weight
%and the weights which occur lie in the rootsystem $R$ of $g$
\[\label{grading}
\qlie{g}=\bigoplus_{\a\in R}{\cal L}_{\a}\oplus{\cal L}_0,~~~~
h\circ a_\a=\a(h) a_\a~\forall a_\a\in{\cal L}_{\a}.
\]
where the dimension of  ${\cal L}_0$ is equal to the rank of $g$
and the ${\cal L}_\a$ are one-dimensional for any root $\a$ of $g$.
\Eq{grading} defines a grading of the quantum Lie algebra:
$[{\cal L}_{\a}\circ{\cal L}_{\b}]\in{\cal L}_{\a+\b}$.
We will refer to
${\cal L}_0={\cal H}$ as the Cartan subalgebra and to the elements
of ${\cal L}_\a$ as root vectors.

We choose some basis $\{X_\a|\a\in\R\}\cup\{H_i|i=1\dots
\mbox{rank(g)}\}$ for the quantum Lie algebra $\qlie{g}$
so that $X_\a\in{\cal L}_\a,
H_i\in{\cal H}$. Because of the grading \Eq{grading}
the Lie bracket relations of $\qlie{g}$ are restricted to
take the form
\[\label{qliestruc}
&&\wb{H_i,X_\a}=l_\a(H_i)\,X_\a,~~~\wb{X_\a,H_i}=-r_\a(H_i)\,X_\a,\non
&&\wb{H_i,H_j}=f_{ij}{}^k\,H_k,~~~~~
\wb{X_\a,X_{-\a}}=-H_\a\in{\cal L}_{0},\\
&&\wb{X_\a,X_\b}=N_{\a\b}\,X_{\a+\b}~~~
\mbox{for }\b\neq -\a \mbox{ and }\a+\b\in R.\nn
\]
This is similar in form to the classical relations \Eq{liestruc}.
There are however some crucial differences.
Because the quantum Lie bracket is not anti-symmetric, there are
two sets of roots, the `left' roots $l_\a$ and the `right' roots
$r_\a$. Furthermore these roots are now not valued in $\CC$ but in
$\Ch $. Similarly the constants $N_{\a,\b}$ and $f_{ij}{}^k$ are elements
of $\Ch $
\footnote{It will usually be possible to treat $\hhh$ as a numeric
deformation parameter
and then to work over $\CC$. However before doing this it is clearly
necessary to verify that all occuring power series in $\hhh$
converge for a certain range of values for $\hhh$.}.
Note also that $\wb{H_i,H_j}$ can be non-zero.

The requirement of invariance of $\qlie{g}$ under $\t{\th}$, $\t{S}$
and $\tau$
is not empty. Already the example of $g=a_2$, treated in section
\ref{sectsl3}, exhibits a whole family of $\uqg$ \mbox{$\circ$ -} submodules
which satisfy the first part of definition \ref{defqlie} but not the
second. However, given any non-invariant
$\uqg$ $\circ$ - submodule $\qlie{g}_0 \subset \uqg$ satisfying the first
part of the definition, a symmetrization with respect to $\t{\th}$,
$\t{S}$ and $\tau$ is always possible. To see this, fix a highest
weight state $\psi \in \qlie{g}_0$ and choose a lowest weight state
$\bar \psi = P(x^-) \circ \psi$, $P(x^-)$ being a monomial in the generators
$x^-_i$. $\t{\th}(\qlie{g}_0)$ satisfies
the first part of the definition as well. Fix a highest weight state
$\psi^{\prime} \in \t{\th}(\qlie{g}_0)$ by requiring
$\t{\th}(\psi^{\prime}) = \bar \psi$ and define
$\bar \psi^{\prime} =  P(x^-) \circ \psi^{\prime}$. Then
$\t{\th}(\psi) = \a \bar \psi^{\prime}$ with some $0 \not= \a \in \Ch$.
The equalities
$\psi = \t{\a} \t{P}(x^+) P(x^-) \circ \psi$ and
$\psi^{\prime} = \a \t{P}(x^+) P(x^-) \circ \psi^{\prime}$
imply $\a = \t{\a}$. Due to the classical limit $\a$ has a square root.
Renormalizing $Q(x^-) = \sqrt{\a} P(x^-)$, $\phi = \psi$,
$\bar \phi = Q(x^-) \circ \phi$, $\phi^{\prime} = \sqrt{\a} \psi^{\prime}$
and $\bar \phi^{\prime} = Q(x^-) \circ \phi^{\prime}$ leads to
$\t{\th}(\phi) = \bar \phi^{\prime}$ and $\t{\th}(\phi^{\prime}) = \bar \phi$.
$\phi_1 = \g \phi + \t{\g} \phi^{\prime}$ satisfies
$\t{\th}(\phi_1) = Q(x^-) \circ \phi_1$ for $0 \not= \g \in \Ch$ arbitrary
but fixed.

Note that the above construction goes through under the additional
requirement $\t{Q}(x^-) = q^{-4 \l(h_\rho)} Q(x^-)$, $\l$ the highest
root. This means $S^2(Q(x^-)) = \t{Q}(x^-)$.
Now let $\phi_2 = \phi_1 - \t{S}(\phi_1)$ to find $\t{S}(\phi_2) = -\phi_2$
and, using (\ref{swb}), $\t{\th}(\phi_2) = Q(x^-) \circ \phi_2$ as
desired. Hence, $\qlie{g} = \uqg \circ \phi_2$ is a quantum Lie algebra.
In case there is a diagram automorphism $\tau$ it is possible
to additionally symmetrise with respect to $\tau$. $\tau$ then
restricts to an automorphism of the resulting quantum Lie algebra.

Whenever there exists more than one quantum Lie algebras associated
to the same Lie algebra $g$, then there exist also whole families of
``almost'' quantum Lie algebras which satisfy the first part of the
definition but are not invariant under $\t{\th}$. Consider the
situation of two quantum Lie algebras with highest weight states
$\psi_1$, $\psi_2$ such that $\t{\th}(\psi_j) = Q(x^-) \circ \psi_j$,
$j=1,2$. (An example for this situation is provided by $a_2$.) For
$\a,\b \in \Ch$ construct the orbit $\uqg \circ (\a \psi_1 + \b \psi_2)$.
{}From $\t{\th}(\a \psi_1 + \b \psi_2)
= Q(x^-) \circ (\t{\a} \psi_1 + \t{\b} \psi_2)$ it follows that
$\uqg \circ (\a \psi_1 + \b \psi_2)$ is $\t{\th}$-invariant only if
$\frac{\a}{\b} = \widetilde{(\frac{\a}{\b})}$.

It had been observed already in the context of the bicovariant differential
calculus that quantum Lie algebras are not left invariant by the
antipode, see \ct{Sch94} for a discussion. We have just shown that
it is however always possible  to find quantum Lie algebras invariant
under the combined action of the antipode and q-conjugation. This
invariance will be crucial in the developments to follow.

There always exist a quantum Lie algebra $\qlie{g}$ associated to any simple 
complex Lie algebra $g$. Furthermore, all $\qlie{g}$ associated to the
same $g$ are isomorphic as algebras. This has been shown in \ct{Del95b}

\subsection{Killing form\lab{sectkilling}}

The Killing form plays a crucial role in the structure theory of
Lie algebras. It is a symmetric bilinear form on the Lie algebra and its
crucial property is the invariance under the adjoint action.
We define a quantum analogue.
\begin{definition}\lab{kill}
The \underline{quantum Killing form} is the map
$B: \dlie\ot\dlie\rightarrow\Ch$
given by
\[\label{qkill}
\kill{a,b}=-Tr_{adj}\left(\t{S}({a})\,b\,u\right).
\]
Here $Tr_{adj}$ denotes the trace over the adjoint representation and
$u$ is the element of $\uqg$ expressing the square of the antipode as
in \Eq{defu}.
\end{definition}
This definition goes over into the classical Killing form $b$ in the
classical limit ($\hhh=0$). From the non-degeneracy of the classical Killing
form the non-degeneracy of the quantum Killing form follows.
The analog of the ad-invariance for the quantum Killing form is
\[\label{adkill}
\kill{a,c\circ b}=\kill{\t{S}(c)\circ a,b}
\]
which can be straightforwardly derived from the definition.

Note that our quantum Killing form on $\qlie{g}$ is not the restriction
of the usual Killing form $R$ on $\uqg$ first defined by Rosso \ct{Ros90}.
The ad-invariance of the later is expressed in terms of the
$\uqg$ coproduct: $\sum\,R(x_{(1)}\circ y,x_{(2)}\circ z)=
\epsilon(x)\,R(y,z)$
 $\forall x,y,z\in\uqg$. This is not useful for our purposes because
the $\uqg$ coproduct leads out of the quantum Lie algebra $\qlie{g}$.

The quantum Killing form is $q$-linear in its first argument and
linear in the second, i.e., for any $\l\in\Ch$
\[\label{qsesqui}
\kill{\l\,b,a}=\t{\l}\,\kill{b,a},~~~~
\kill{b,\l\,a}=\l\,\kill{b,a}.
\]
The quantum Killing form is not symmetric. However it is $q$-symmetric in
the sense that
\[\label{ks}
\kill{b,a}=\widetilde{\kill{a,b}}
\]
In addition we have the following two relations
\[
\label{kss}
\kill{b,a}&=&\kill{\t{S}(a),S(\t{b})}\\
\label{ks1}&=&\kill{\t{\th}(a),\t{\th}(b)}.
\]
To derive these relations one has to realize that the dual of
$\pi_{adj}$, $\pi_{adj}\cdot\th$ and $\sim\cdot\pi_{adj}\cdot\sim$
are all related to the adjoint representation $\pi_{adj}$ itself by
similarity transformations.

\subsection{Weyl canonical form\lab{sectcanon}}

\begin{proposition}\lab{canon}
It is possible to choose a basis $\{X_\a|\a\in\R\}\cup\{H_i|i=1\dots
\mbox{rank(g)}\}$ for $\dlie$ with the properties
\[
\label{kxx}&&\kill{{X}_\a,X_{-\a}}=-1\\
\label{ktx}&&\t{\th}(X_\a)=X_{-\a},~~~~\t{\th}(H_i)=-H_i,\\
\label{ksx}&&\t{S}(X_\a)=-q^{-\cp{\rho,\a}}\,X_\a,~~~~\t{S}(H_i)=-H_i.
\]
\end{proposition}

{\it Proof}. In the following we will have to invert and to take square
roots of elements of $\Ch$. While this is not in general possible, it
presents no problem for those formal power series which have a
nonvanishing classical limit.

It is clear by (\ref{adkill}) that $\kill{{X}_\a,X_\b}\propto\d_{\a+\b,0}$.
We can choose the scale of the the $X_\a$ for positive $\a$ so
that $\kill{{X}_{-\a},X_{\a}}=-1$. By the symmetry property \Eq{ks}
of the Killing form the normalization \Eq{kxx} then holds for all $\a$.

The most general action of $\th$ on $X_\a$ is, for reasons of weight,
$\th(X_\a)=f_\a\,\t{X}_{-\a}$ for some $f_\a\in\CC[[\hhh]]$. Since
$\th^2=\mbox{id}$ and $\sim$ commutes with $\th$ we have
$X_\a = f_\a \t{f}_{-\a} X_\a$, i.e. $f_\a^{-1}=\t{f}_{-\a}$
for all $\a$. If we rescale the $X_\a$ by $f_\a^{-1/2}$ both
\Eq{ktx} and \Eq{kxx} hold. In particular, the $X_\a$ are determined up to
sign by \Eq{kxx} and \Eq{ktx}.

The most general action of $S$ on $X_\a$ is, again for reasons of weight,
$S(X_\a)=s_\a\,\t{X}_\a$ for some $s_\a\in\CC[[\hhh]]$. Hence
$X_\a = S^{-1} \cdot S(X_\a) =
s_\a \widetilde{S(X_\a)} = s_\a \t{s}_\a X_\a$, i.e.
$s_\a^{-1}=\t{s}_\a$. Furthermore, $S\cdot\th\cdot S=\th$ and
\Eq{ktx} imply $s_{-\a}=s_\a^{-1}$.
Finally, \Eq{defu}, \Eq{qsesqui} and \Eq{kss} lead to
$1=-\kill{{X}_{-\a},X_\a} = -\kill{\t{S}(X_\a),S(\t{X}_{-\a})} =
s_\a\,\t{s}_{-\a}\, q^{-2\cp{\rho,\a}}$. Hence $s_\a^2=q^{2\cp{\rho,\a}}$.
The sign in \Eq{ksx} is determined by the classical limit ($\hhh=0$).

We construct the basis states for the Cartan subalgebra ${\cal H}$ as follows
\[\label{hbasis}
H_i&=&\half\left(q^{\rho\cdot\a_i}\wb{X_{-\a_i},X_{\a_i}}-
q^{-\rho\cdot\a_i}\wb{X_{\a_i},X_{-\a_i}}\right).
\]
Then, using \Eq{bwr} and \Eq{tiso},
$\th(H_i)=-\t{H}_i$ follows.
The relation $S(H_i)=-\t{H}_i$ follows from
$\kill{X_{-\a},\wb{H_i,X_\a}} =
\kill{\wb{\t{\th}(H_i),\t{\th}(X_\a)},\t{\th}(X_{-\a})} =
- \kill{\wb{H_i,X_{-\a}},X_\a} =
- \kill{X_{-\a},\wb{\t{S}(H_i),X_\a}}$.
At $\hhh=0$ the $H_i$ defined above are equal to the
$\hat{h}_i$ of the usual Weyl canonical form of \Eq{liestruc}.
This shows that the $H_i$ are linearly independent and thus give a
basis of the Cartan subalgebra ${\cal H}.$ $\Box$

{\it Remark}. Note that the $H_i$ are not unique. For example, every
choice
$H_i = \half\left(\g_i H_{-\a_i} - \t{\g_i} H_{\a_i}\right)$
with $0 \not= \g_i \in \Ch$ such that
$\g_i + \hhh \Ch = 1 + \hhh \Ch$ is possible.

If $g$ has a diagram automorphism then $\tau$
acts as
\[\label{tau}
\tau(X_\a)=t_\a\, X_{\tau(\a)},~~~~~~~
\tau(H_i)=H_{\tau(i)},~~~~~~t_\a=\pm 1,
\]
where the signs $t_\a$ are the same as in the classical case.
{\it Proof}. For reasons of weight $\tau(X_\a)=t_\a\,X_{\tau(\a)}$ for some
$t_\a\in\Ch$. From $\kill{X_{-\a},X_\a}=\kill{\tau(X_{-\a}),\tau(X_\a)}=
\t{t}_{-\a}t_\a\kill{X_{-\a},X_\a}$ it follows that $\t{t}_{-\a}t_\a=1$.
{}From $\t{\th}(\tau(X_{-\a}))=\tau(\t{\th}(X_{-\a}))$ it follows that
$\t{t}_{-\a}=t_\a$. Together this gives $t_\a^2=1$ and thus
$t_\a=\pm 1$. The action on $H_i$ follows from \Eq{hbasis} and the
choice $t_{\a_i}=1.~\Box$

\subsection{Relations between structure constants\lab{sectrelations}}

We are now ready to derive relations between the various
structure constants appearing in \Eq{qliestruc}
when using the basis of proposition \ref{canon}.
{}From the isomorphism property \Eq{tiso} of $\t{\th}$ we obtain
an expression of the quantum roots for negative $\a$ in terms of the
quantum roots for positive $\a$
\[
\label{llt}\wb{\t{\th}(H_i),\t{\th}(X_{-\a})}=\t{\th}(\wb{H_i,X_{-\a}})
&~\Rightarrow~&l_{-\a}(H_i)=-\widetilde{l_\a({H}_i)}~~\forall\a,i,\\
\label{rrt}\wb{\t{\th}(X_{-\a}),\t{\th}(H_i)}=\t{\th}(\wb{X_{-\a},H_i})
&~\Rightarrow~&r_{-\a}(H_i)=-\widetilde{r_\a({H}_i)}~~\forall\a,i.
\]
Thus, unlike in the classical case, the negative of a left quantum root is
not a left quantum root again, but the q-conjugated negative is. Idem for
right quantum roots.
We also obtain relations for the structure constants $N$ and $f$
\[
\label{n1}\wb{\t{\th}(X_\a),\t{\th}(X_\b)}=\t{\th}(\wb{X_\a,X_\b})
&~\Rightarrow~&N_{\a,\b}=\t{N}_{-\a,-\b}~~~\forall \a,\b,\\
\label{f1}\wb{\t{\th}(H_i),\t{\th}(H_j)}=\t{\th}(\wb{H_i,H_j})
&~\Rightarrow~&f_{ij}{}^k=-\t{f}_{ij}{}^k~~~\forall i,j,k.
\]
{}From the ad-invariance \Eq{adkill} of the quantum Killing form
we obtain the characterization
of the Cartan subalgebra elements $H_\a$ in terms of the right roots
\[\label{rbh}
-\kill{H_\a,H}&=&\kill{\wb{X_\a,X_{-\a}},H}
=\kill{X_{-\a},\wb{\t{S}(X_\a),H}}\non
&=&\kill{X_{-\a},-q^{-\rho\cdot\a}\wb{X_\a,H}}
=q^{-\cp{\rho,\a}}r_\a(H)~~~~~
\forall\a,\forall H\in{\cal H}.
\]
Because of the non-degeneracy of the Killing form these relations
determine the $H_\a$ uniquely in terms of the roots. We also obtain
further relations for the structure constants $N$ and $f$
\[
\label{n2}&&\kill{\wb{\t{S}(X_\a),{X}_\b},X_{-\a-\b}}=
\kill{{X}_\b,\wb{X_\a,X_{-\a-\b}}}
\non
&&~~~~~~~~~~~~~~\Rightarrow~N_{\a,-\a-\b}=
-q^{\cp{\rho,\a}}\t{N}_{\a,\b}~~~
\forall\a,\b,\\
\label{f2}&&\kill{\wb{\t{S}(H_j),{H}_i},H_k}=
\kill{{H}_i,\wb{H_j,H_k}}
\non
&&~~~~~~~~~~~~~~\Rightarrow~\sum_l f_{jk}{}^lB_{il}=
-\sum_l \t{f}_{ji}{}^l B_{lk},
\]
where we have defined $B_{ij}=\kill{{H}_i,H_j}$.

There exists a quantum Lie algebra anti-automorphism 
$\chi:\qlie{G}\rightarrow\qlie{g}$
acting on the basis as
\[\label{chi}
\chi(X_\a)=-X_{-\a},~~~~~~\chi(H_i)=H_i.
\]
{}From the anti-automorphism property
\[\label{chiauto}
\wb{\chi(a),\chi(b)}=\chi(\wb{b,a})~~~~~~\forall a,b\in\dlie,
\]
we obtain the relation between the `left' and `right' quantum roots
\[\label{leftright}
l_\a=-r_{-\a}~~~~~~\forall\a,
\]
and the relations
\[\label{nf3}
N_{\a,\b}=-N_{-\b,-\a},~~~~~~~~f_{ij}{}^k=f_{ji}{}^k.
\]
The proof that \Eq{chi} defines an anti-automorphism of the
quantum Lie algebras is contained in \ct{Del95b}.

If $g$ has a diagram automorphism $\tau$ then this leads to further
relations
\[\label{diarelsf}
&&f_{\tau(i)\tau(j)}{}^{\tau(k)}=f_{ij}{}^k,~~~~~~
N_{\tau(\a),\tau(\b)}=t_\a t_\b t_{\a+\b} N_{\a,\b},\\
\label{diarelsr}&&l_{\tau(\a)}(H_{\tau(i)})=l_\a(H_i),~~~~~~
r_{\tau(\a)}(H_{\tau(i)})=r_\a(H_i),\\
\label{diarelsb}&&B_{\tau(i)\tau(j)}=B_{ij}\nonumber.
\]

\subsection{Quantum root spaces\lab{sectroots}}

We have seen that a quantum Lie algebra posesses two sets of quantum
roots, $l_\a$ and $r_\a$, defined by
\[
\wb{H,X_\a}=l_\a(H)\,X_\a,~~~~~
\wb{X_\a,H}=-r_\a(H)\,X_\a.
\]
The roots are linear forms on the Cartan subalgebra ${\cal H}$
with values in $\Ch$, i.e.,
they are elements of ${\cal H}^*$.
If the quantum Lie algebra has the anti-automorphism $\chi$
of \Eq{chi}, then the roots are related by $r_\a=-l_{-\a}$, i.e.,
the set of right roots is just the negative of the set of left roots.
%The anti-automorphism exists in the examples we have studied and
%one might expect it to exist in all quantum Lie algebras.

{}From the Killing form on ${\cal H}$ we construct a form on ${\cal H}^*$
in the usual way. To any element $v\in{\cal H}^*$ we associate the
unique element $H_v\in{\cal H}$ satisfying $v(H)=B(H_v,H)~
\forall H\in{\cal H}$. Note that this pairing is $q$-linear in the
sense that the element of {\cal H}
associated to $\l\,v$ for some $\l\in\Ch$ is
not $\l\,H_v$ but $\t{\l}\,H_v$. The form on ${\cal H}^*$ is defined by
\[\label{qform}
\form{v,w}=\kill{H_v,H_w}~~~~~\forall v,w\in{\cal H}^*
\]
Because the Killing form is $q$-linear in the first factor and linear
in the second, the form $\form{\,.\,,\,.\,}$ is linear in the first
factor and $q$-linear in the second
\[
\form{\l\,v,w}=\l\form{v,w},~~~~~
\form{v,\l\,w}=\t{\l}\form{v,w}.
\]
It is also $q$-symmetric
\[
\form{v,w}=\widetilde{\form{w,v}}
\]
{}From the relation \Eq{rbh} we can read off that, for example,
\[
\form{r_\a,r_\b}=q^{\rho\cdot(\a-\b)}\,\kill{H_\a,H_\b}.
\]

In the classical case of complex Lie algebras one introduces a
real form ${\cal H}_{\RR}$ of the Cartan subalgebra and on its dual
space ${\cal H}_{\RR}^*$, which is a real vectorspace, the form
induced by the Killing form is a real, positive definite, bilinear
form, thus giving ${\cal H}_{\RR}^*$ the structure of a Euclidean
space. This is the root space.

We can imitate this construction for quantum Lie algebras.
We define the `q-real' form ${\cal H}_{\RR[[\hhh^2]]}$ of the
Cartan subalgebra as the module over $\RR[[\hhh^2]]$ spanned by
the $H_i$. We choose $\RR[[\hhh^2]]$ as the base ring because it
consists of the elements of $\Ch$ which are invariant under both
complex conjugation and $q$-conjugation. The roots, when restricted
to ${\cal H}_{\RR[[\hhh^2]]}$ still give values in $\RR[[\hhh]]$,
and not in $\RR[[\hhh^2]]$, and thus do not lie in
$\left({\cal H}_{\RR[[\hhh^2]]}\right)^*$. The $q$-symmetrized combinations
$a_\a=\half(r_\a-r_{-\a})$ do however give values in $\RR[[\hhh^2]]$.
The $a_i\equiv a_{\a_i}$ for all simple roots $\a_i$ form a basis for
${\cal H}^*_{\RR[[\hhh^2]]}$. On this basis the form is given by
\[
\form{a_i,a_j}=\kill{H_i,H_j}=B_{ij}.
\]
We see immediately that
the form $\form{\,.\,,\,.\,}$ restricted to ${\cal H}^*_{\RR[[\hhh^2]]}$
is a symmetric, non-degenerate, bilinear form with values in
$\RR[[\hhh^2]]$.

We expect however that in an eventual axiomatic description of quantum
root systems the unrestricted form $\form{\,.\,,\,.\,}$ will be used
and that the fact that it is not symmetric and bilinear but rather
q-symmetric and q-bilinear will play a central role.

\section{Explicit examples\lab{sectex}}

We have explicitly constructed three examples of quantum Lie algebras,
namely those associated to $g=a_2, a_3, c_2$ and to $g_2$. The construction
follows straigthforwardly from the definition \ref{defqlie}.
We search for a highest weight state inside $U_\hhh^{\ge 0} (g)$ and impose
a symmetry constraint if appropriate. Then the corresponding orbit is
constructed and explicitly tested for the invariance properties required by
the definition to be satisfied. The details for the cases of $g=a_2$ and $c_2$
are given below.

Rather than to describe the quantum groups in terms of fundamental generators
and  their relations, the selection of a Poincare Birkoff Witt (PBW) type basis
is useful for explicit computations. In the construction of
such a basis with the help of the Lusztig automorphisms
\ct{Lus93} we follow the conventions
of \ct{Cha94}; an alternative would be the approach of \ct{Kho91}. For a
reduced decomposition of the longest Weyl group element $w_0 =
s_{i_1} \dots s_{i_N}$ the quantum root vectors are given by
\[
e_k &=& T_{i_1} \dots T_{i_{(k-1)}} (X_{i_{k}}^+) \quad \mbox{and} \\
f_k &=& T_{i_1} \dots T_{i_{(k-1)}} (X_{i_{k}}^-).
\]
Note that $e_k$ is a polynomial in $\lbrace X_i^+ \rbrace$ while analogously
$f_k$ is a polynomial in $\lbrace X_i^- \rbrace$, although this is not entirely
obvious from the definition of the Lusztig automorphisms $T_j$.

\subsection{$\qlie{a_2}$\lab{sectsl3}}

$a_2=sl_3$ is the rank 2 Lie algebra with Cartan matrix
\[
a=\left(\begin{array}{cc}2&-1\\-1&2\end{array}\right)
\]
It has a diagram automorphism $\tau$ which exchanges the two
simple roots, i.e.,
$X_1^\pm \leftrightarrow X_2^\pm$ and $k_1 \leftrightarrow k_2$.
The quantum root vectors generating the PBW basis which we use involve the
choice $w_0 = s_1 s_2 s_1$:
\[
e_1&=&X_1^+,~~~~e_2=-X_1^+X_2^++q^{-1}\,X_2^+X_1^+,~~~~
e_3=X_2^+,\non
f_1&=&X_1^-,~~~~f_2=q\,X_1^-X_2^- -X_2^-X_1^-,~~~~~~~
f_3=X_2^-.
\]
The diagram automorphism acts as
$\tau(e_1) = e_3$,
$\tau(e_2) = -q^{-1} e_2 - (1-q^{-2}) e_3 e_1$,
$\tau(e_3) = e_1$.
In terms of the PBW basis it is straightforward to write down
an Ansatz $\Psi$ for a highest weight state according to point 2) of the
definition \ref{defqlie}. Once we restrict the Ansatz for $\Psi$ to lie
entirely in $U^{\geq 0}_\hhh$, i.e.,not to contain any $f_i$, we find two
independent solutions of the equations $x_i^+ \circ \Psi = 0$.
With respect to the diagram automorphism these can be described as a highest
weight state
\[
\Psi^+ = e_2 (k_1^{1/3} k_2^{-1/3} - q^{-1} k_1^{-1/3} k_2^{1/3})
- (1 - q^{-2}) e_3 e_1 k_1^{-1/3} k_2^{1/3}
\]
that is invariant under the diagram symmetry, while
\[
\Psi^- = e_2 (k_1^{1/3} k_2^{-1/3} + q^{-1} k_1^{-1/3} k_2^{1/3})
+ (1 - q^{-2}) e_3 e_1 k_1^{-1/3} k_2^{1/3}
\]
changes sign under the diagram automorphism.  The (skew)invariance of
$\Psi^\pm$ follows by means of $\left[e_1,e_3\right]_{q^{-1}} = -e_2$.
The symmetrisation with respect to $\tau$ enforces the symmetries
required by the definition of a quantum Lie algebra.

We now observe that $\Psi^+$ vanishes in the classical limit $q \to 1$ whilst
$\Psi^-$ reduces to the highest root vector of the classical Lie algebra. Hence
$\Psi^-$ is a desirable starting point for the construction of an adjoint
orbit. The resulting orbit is in fact found to satisfy all the
requirements of definition \ref{defqlie}. We then chose a quantum
Weyl basis with the properties of proposition \ref{canon}. The explicit
expressions for these quantum Lie algebra generators are listed below
to give the reader an idea about the form of these generators. Note for
example that the quantum Cartan subalgebra generators are not simple
expressions.
\[\label{gensl3}
X_{\a_1+\a_2} &=& -C\left(e_2 (q^{-1/2} k_1^{1/3} k_2^{-1/3} +
q^{-3/2} k_1^{-1/3} k_2^{1/3})
- (q^{-1}-q)q^{-3/2} e_3 e_1 k_1^{-1/3} k_2^{1/3}\right) \non
X_{\a_2} &=& -iC\left( e_3 (q^{1/2}k_1^{-2/3} k_2^{-1/3} +
q^{-1/2} k_1^{2/3} k_2^{1/3})
+ (q^{-1}-q)q^{-1/2} e_2 f_1 k_1^{1/3} k_2^{-1/3}\right) \non
X_{\a_1} &=&iC\left( e_1 (q^{1/2}k_1^{-1/3} k_2^{-2/3} +
 q^{-1/2} k_1^{1/3} k_2^{2/3})
-(q^{-1}-q)q^{-3/2} e_2 f_3 k_1^{-1/3} k_2^{1/3}\right. \non
&& \mbox{} \hspace{70pt} \left.
+ (q^{-1}-q)^2 q^{-3/2} e_3 e_1 f_3 k_1^{-1/3} k_2^{1/3}\right) \non
H_1 &=&C^2 \frac{1+q^3}{2(1-q)}
\left(-q k_1^{2/3} k_2^{-2/3} + k_1^{-2/3} k_2^{+2/3}
- k_1^{4/3} k_2^{2/3} + q k_1^{-4/3} k_2^{-2/3}\right.
\non
&&~~~~~~~~~~~~~\left.+(1-q^2)^2\left(
+ e_1 f_1 (q^{-2}k_1^{-1/3} k_2^{-2/3} + q^{-3} k_1^{1/3} k_2^{2/3})
\right.\right.\non
&&\left.\left.~~~~~~~~~~~~~~~~~~~~~~~~~~~~
+q^{-4} e_2 f_2 k_1^{-1/3} k_2^{1/3}
-q^{-2} e_3 f_3 k_1^{2/3} k_2^{1/3}\right)\right.
\non
&&~~~~~~~~~~~~~\left.-(1-q^2)^3
 q^{-5} e_3 e_1 f_2 k_1^{-1/3} k_2^{1/3} \right) \non
H_2 &=&C^2 \frac{1+q^3}{2(1-q)}
\left(-q k_1^{-2/3} k_2^{2/3} + k_1^{2/3} k_2^{-2/3}
- k_1^{2/3} k_2^{4/3} + q k_1^{-2/3} k_2^{-4/3}\right.
\non
&&~~~~~~~~~~~~~\left.+(1-q^2)^2\left(
+ e_3 f_3 (q^{-2}k_1^{-2/3} k_2^{-1/3} + q^{-3} k_1^{2/3} k_2^{1/3})
\right.\right.\non
&&\left.\left.~~~~~~~~~~~~~~~~~~~~~~~~~~~~
+q^{-4} e_2 f_2 k_1^{1/3} k_2^{-1/3}
-q^{-2} e_1 f_1 k_1^{1/3} k_2^{2/3}\right)\right.
\non
&&~~~~~~~~~~~~~\left.+(1-q^2)^3
 q^{-4} e_2 f_3 f_1 k_1^{1/3} k_2^{-1/3} \right) \non
X_{-\a_1} &=& iC\left( f_1 (q^{1/2}k_1^{2/3} k_2^{-2/3} +
q^{-1/2} k_1^{4/3} k_2^{2/3})
- (q^{-1}-q)q^{-1/2} e_3 f_2 k_1^{2/3} k_2^{1/3}\right) \non
X_{-\a_2} &=& -iC\left( f_3 (q^{1/2}k_1^{-2/3} k_2^{2/3} +
 q^{-1/2} k_1^{2/3} k_2^{4/3})
+(q^{-1}-q)q^{-3/2} e_1 f_2 k_1^{1/3} k_2^{2/3}\right. \non
&& \mbox{} \hspace{70pt} \left.
+ (q^{-1}-q)^2 q^{-1/2} e_1 f_3 f_1 k_1^{1/3} k_2^{2/3}\right) \non
X_{-\a_1-\a_2} &=& C\left(f_2 (q^{-1/2} k_1^{2/3} k_2^{4/3} +
q^{-3/2} k_1^{4/3} k_2^{2/3})
+(q^{-1}-q)q^{-1/2} f_3 f_1 k_1^{4/3} k_2^{2/3}\right)\nonumber
\]
The normalisation factor is
\[
C=\left(2(q^{-1/2}+q^{1/2})(q^{-3/2}+q^{3/2})
(q^{-3}+q^{-1}-1+q+q^3)\right)^{-1/2}.
\]
It could be absorbed into a different normalization of the quantum Killing
form in \Eq{qkill}.

The left quantum roots are, using the notation $H=\sum\,h_i\,H_i$,
\[\label{lsl3}
l_{\a_1}(H)&=&l \left((q^{-3/2}+q^{-1/2})\,h_1
-q^{1/2}\,h_2\right)\non
l_{\a_2}(H)&=&l \left(-q^{1/2}\,h_1+(q^{-3/2}+q^{-1/2})\,h_2
\right)\non
l_{\a_1+\a_2}(H)&=&l \,q^{-3/2} \,(h_1+h_2)\non\\
&&~~~~~~l=C^2 (q^{-1/2}+q^{1/2})(q^{-3/2}+q^{3/2})^2/2\nn
\]
The negative roots are obtained by $q$-conjugation according to
\Eq{llt}. The right roots are given according to
\Eq{leftright}. The roots are seen to be related by the
diagram automorphism according to \Eq{diarelsr}

The q-conjugation-invariant roots
$a_\a=\frac{1}{2}(r_\a + l_\a)$
introduced in section \ref{sectroots} are
\[\label{asl3}
a_{\a_1}&=&\frac{l}{2}\left((q^{-3/2}+q^{-1/2}+q^{1/2}+q^{3/2})\,h_1
-(q^{-1/2}+q^{1/2})\,h_2\right)\non
a_{\a_2}&=&\frac{l}{2}\left(-(q^{-1/2}+q^{1/2})\,h_1
+(q^{-3/2}+q^{-1/2}+q^{1/2}+q^{3/2})\,h_2\right)\non
a_{\a_1+\a_2}&=&\frac{l}{2}\,(q^{-3/2}+q^{3/2})\,(h_1+h_2)
\]
These have the classical properties
\[\label{latsl3}
a_\a+a_\b=a_{\a+\b},~~~~~~a_{-\a}=-a_\a,
\]
i.e., they form a root lattice. This interesting feature, which
makes these root systems look very similar to their classical
counterparts, is true for $g=a_n$ for any $n$ \ct{Del95}, but is
not true for $c_2$, as we will see in the next section.

The Killing form on the Cartan subalgebra is given by the matrix
$B$ with entries $B_{ij}=\kill{H_i,H_j}=\form{a_{\a_i},a_{\a_j}}$,
\[\label{killsl3}
B&=&b\,\left(
\begin{array}{cc}
q+q^{-1}&-1\\
-1&q+q^{-1}
\end{array}\right)\non
b&=&\left((q^{-1/2}+q^{1/2})^2(q^{-3/2}+q^{3/2})^2\right)\,C^2/4
\]
The pairwise equality of the elements is due to the
diagram automorphism.

Once one has the knowledge of the Killing form and of the
roots,
the $H_\a$, which appear as the result of $\wb{X_\a,X_{-\a}}$, are
determined by \Eq{rbh}. In terms of the $H_i$ they read
\[
H_{\a_1}&=&a\left(-q^{-1/2}\,H_1+(-q^{1/2}+q^{3/2})H_2\right)\non
H_{\a_2}&=&a\left((-q^{1/2}+q^{3/2})H_1-q^{-1/2}\,H_2\right)\non
H_{\a_1+\a_2}&=&-a\,q^{1/2}(H_1+H_2)\non
H_{-\a_1}&=&a\left(q^{1/2}\,H_1+(q^{-1/2}-q^{-3/2})H_2\right)\non
H_{-\a_2}&=&a\left((q^{-1/2}-q^{-3/2})H_1+q^{1/2}\,H_2\right)\non
H_{-\a_1-\a_2}&=&a\,q^{-1/2}(H_1+H_2)\non
&&~~~~~~a=2(q^{-3/2}+q^{3/2})^{-1}
\]
Note that the coefficients in the expansion of the $H_\a$ are related
to those in $H_{-\a}$ by q-conjugation and sign change.

We need to give only one of the structure constants $N$
\[
N_{\a_1,\a_2}=(q^{-3/2}+q^{3/2})\,C
\]
Through the relations \Eq{n2} and \Eq{nf3} all the other non-zero $N_{\a,\b}$
are related to this (note that $N_{\a_1,\a_2}=\t{N}_{\a_1,\a_2}$)
\[
&&
N_{\a_1,-\a_1-\a_2}=-q\,N_{\a_1,\a_2},~~
N_{\a_2,\a_1}=-\,N_{\a_1,\a_2},~~
N_{\a_2,-\a_1-\a_2}=q\,N_{\a_1,\a_2},~~
\non
&&
N_{\a_1+\a_2,-\a_1}=q\,N_{\a_1,\a_2},~~
N_{\a_1+\a_2,-\a_2}=-q\,N_{\a_1,\a_2},~~
N_{-\a_1,-\a_2}=\,N_{\a_1,\a_2},~~
\non
&&
N_{-\a_1,\a_1+\a_2}=-q^{-1}\,N_{\a_1,\a_2},~~
N_{-\a_2,\a_1+\a_2}=q^{-1}\,N_{\a_1,\a_2},~~
N_{-\a_2,-\a_1}=-\,N_{\a_1,\a_2}.
\non
&&N_{-\a_1-\a_2,\a_1}=q^{-1}\,N_{\a_1,\a_2},~~
N_{-\a_1-\a_2,\a_2}=-q^{-1}\,N_{\a_1,\a_2},~~
\]
This is confirmed by the results of the explicit calculations.

For the structure constants $f_{ij}{}^k$ for the Cartan subalgebra we find
\[
&&f_{11}{}^1=f_{22}{}^2=-f\, (q^{-2}+q^{-1}+1+q+q^2)\non
&&f_{22}{}^1=f_{11}{}^2=-f\,(q^{-1}+q)\non
&&f_{12}{}^1=f_{21}{}^1=f_{12}{}^2=f_{21}{}^2=f\non
&&~~~~~~~~~f=(q^{1/2}-q^{-1/2})(q^{-1/2}+q^{1/2})^2(q^{-3/2}+q^{3/2})\,C^2/2
\]

\subsection{$\qlie{c_2}$\lab{sectsp4}}

$c_2=sp(4)=b_2=so(5)$ is the rank 2 Lie algebra with Cartan matrix
\[
a=\left(\begin{array}{cc}2&-2\\-1&2\end{array}\right)
\]
It has no diagram automorphisms. With conventions analogous to the
previous example
\[
l_{\a_1}(H)&=&l\left(
(q^{-2}-1+q^2)^2 q^{-1}\,h_1-q^3\,h_2\right)
\non
l_{\a_2}(H)&=&l\left(
-(q^{-2}-1+q^2)q^{-1}\,h_1+(q^{-1}+q)q^{-2}\,h_2\right)
\non
l_{\a_1+\a_2}(H)&=&l\left(
(q^{-1}-q)(q^{-2}-1+q^2)q^{-2}\,h_1+q^{-1}\,h_2\right)
\non
l_{2\a_1+\a_2}(H)&=&l\left(
(q^{-2}-1+q^2)q^{-3}\,h_1\right)
\non\vspace{2mm}
&&l=(q^{-1}+q)^3(q^{-2}-1+q^2)^2\,C^2/2
\]
\[\label{asp4}
a_{\a_1}(H)&=&a\left(
(q^{-2}-1+q^2)^2\,h_1-(q^{-2}-1+q^2)\,h_2\right)
\non
a_{\a_2}(H)&=&a\left(
-(q^{-2}-1+q^2)\,h_1+(q^{-2}+q^2)\,h_2\right)
\non
a_{\a_1+\a_2}(H)&=&a\left(
(q^{-1}-q)^2(q^{-2}-1+q^2)\,h_1+h_2\right)
\non
a_{2\a_1+\a_2}(H)&=&a\left(
(q^{-2}-1+q^2)^2\,h_1\right)
\non\vspace{2mm}
&&a=(q^{-1}+q)\,C^2/2
\]
The normalisation constant is
\[
C&=&\left((q^{-1}+q)^2(q^{-2}+q^2)(q^{-1}+1+q)(q^{-1}-1+q)\right.\non
&&~~~~\left.(q^{-2}-1+q^2)(q^{-4}-q^{-2}+1-q^2+q^4)\right)^{-1/2}.
\]
Note that
\[
a_{2\a_1+\a_2}\neq a_{\a_1+\a_2}+a_{\a_1}
\]
Thus, in contrast to the case of $g=a_n$, these roots do not form
a root lattice.\footnote{It is tempting to speculate that there may be a
relation between the non-closure of the above root triangle and
the non-closure of some mass triangles in the affine Toda theory
based on $c_2$.}

The Killing form on the Cartan subalgebra is given by the matrix $B$ with
entries $B_{ij}=\kill{H_i,H_j}=\form{a_{\a_i},a_{\a_j}}$,
\[
B&=&b\,\left(
\begin{array}{cc}
q^{-2}-1+q^2&-1\\
-1&
\frac{q^{-2}+q^2}{q^{-2}-1+q^2}
\end{array}\right)\non
b&=&\left((q^{-1}+q^{1})^4
(q^{-2}-1+q^2)^3\right)\,C^2/4
\]

We find
\[
N_{\a,\b}=-(q^{-1}+q)(q^{-2}-1+q^2)\,n_{\a,\b}
\]
with the $n_{\a,\b}$ given in the following table. The rows are labeled
by $\a$ and the coloums by $\b$.
\[
\begin{array}{r|cccccccc}
&\begin{array}{c}2\a_1\\~+\a_2\end{array}&
\begin{array}{c}\a_1\\~+\a_2\end{array}&
{}~~\a_2&~~\a_1&~-\a_1&~-\a_2&
\begin{array}{c}-\a_1\\~-\a_2\end{array}&
\begin{array}{c}-2\a_1\\~-\a_2\end{array}\\
\hline
2\a_1+\a_2&0&0&0&0&-q^2&0&q^2&0\\
\a_1+\a_2&0&0&0&q&q^3&-1&0&-q^2\\
\a_2&0&0&0&-q^2&0&0&1&0\\
\a_1&0&-q^{-1}&q^{-2}&0&0&0&-q^3&q^2\\
-\a_1&q^{-2}&-q^{-3}&0&0&0&q^2&-q&0\\
-\a_2&0&1&0&0&-q^{-2}&0&0&0\\
-\a_1-\a_2&-q^{-2}&0&-1&q^{-3}&q^{-1}&0&0&0\\
-2\a_1-\a_2&0&q^{-2}&0&-q^{-2}&0&0&0&0
\end{array}
\]
In view of the relations given in section \ref{sectrelations} the
structure constants are fixed once $N_{\a_1,\a_2}$ and
$N_{\a_1,\a_1+\a_2}$ are given.
Also because of these relations the above table is q-anti-symmetric
about the diagonal and antisymmetric about the opposite diagonal.

For the structure constants $f_{ij}{}^k$ for the Cartan subalgebra we find
\[
f_{11}{}^1&=&-f\,(q^{-2}-1+q^2)(q^{-4}-q^{-2}+3-q^2+q^4)\non
f_{22}{}^2&=&-f\, (q^{-2}-q^{-1}+1-q+q^2)(q^{-2}+q^{-1}+1+q+q^2)\non
f_{11}{}^2&=&-f\,(q^{-2}+q^2)(q^{-2}-1+q^2)^2\non
f_{22}{}^1&=&-f\,(q^{-2}+q^2)(q^{-2}-1+q^2)^{-1}\non
f_{12}{}^2&=&f_{21}{}^2=f\,(q^{-2}+q^2)(q^{-2}-1+q^2)\non
f_{12}{}^1&=&f_{21}{}^1=f\non
\vspace{2mm}
&&f=-(q^{-1}-q)(q^{-1}+q)^3(q^{-2}-1+q)/2
\]

The the quantum roots corresponding to the positive classical
roots are represented by
\[
H_{\a_1}&=&d\left((-q^{-4}+q^{-2}-2+q^2)\,H_1-\right.\non
&&~~~\left.(q^{-1}-q)(q^{-2}+q^2)(q^{-2}-1+q^2)q^{-1}\,H_2\right)
\non
H_{\a_2}&=&d\left(-(q^{-1}-q)q^2\,H_1-
(q^{-2}-1+q^2)q^{-1}\,H_2\right)
\non
H_{\a_1+\a_2}&=&d\left((-q^{-2}+1-2q^2-q^4)\,H_1-
(q^{-2}-1+q^2)\,H_2\right)
\non
H_{2\a_1+\a_2}&=&d\left(-(q^{-2}+q^2)q\,H_1-
(q^{-2}-1+q^2)q\,H_2\right)
\non\vspace{2mm}
&&d=2\left((q^{-1}+q)(q^{-2}-1+q^2)\right)^{-1}
\]
Again the corresponding expressions for negative $\a$ are obtained by
q-conjugating the coefficients and changing the sign.

\section{Discussion \lab{sumdisc}}

We have shown that it is possible to develop a theory of quantum
Lie algebras in terms of an analogue of Weyl's canonical form and the
resulting quantum roots and structure constants. The key idea is the
concept of q-conjugation that allows us to exploit q-linear analogues
of the antipode and the Cartan involution in connection with a generalised
Killing form.

Objects similar to our quantum Lie algebras have been studied in the
framework of bicovariant differential calculus on quantum groups,
see \ct{Asc93} for a very readable review. There one considers the
dual space to the space of left-invariant one-forms, which is a
$\circ$-submodule of $\uqg$\footnote{Rather than working with modules
over $\Ch$ people treat $q=\exp(\hhh)$ as a number and work with
vector spaces over $\CC$ or $\RR$.}. The case $g=sl_3$ has been
explicitly worked out in \ct{Asc92}. It does not coincide with
our quantum Lie algebra $\qlie{sl_3}$ studied in section \ref{sectsl3}.
In particular the module of \ct{Asc92} is not invariant under the
diagram automorphism of $sl_3$. 

%For the case of $g=a_n$, we
%find that his module fails to be a quantum Lie algebra in the sense of
%definition
%\ref{defqlie} only because it is not invariant under $\t{\th}$ and $\t{S}$.
%One could equally well consider the dual space to the space of
%right-invariant one-forms with respect to the opposite coproduct, which
%is also a $\circ$-submodule of $\uqg$. It too fails to give a
%quantum Lie algebra according to definition \ref{defqlie} because it is
%not invariant under $\t{\th}$ and $\t{S}$. Furthermore, these two
%``almost'' quantum Lie algebras are not isomorphic to each other.
%Our requirement of invariance under the action of the tilded antipode
%$\t{S}$ selects special linear combinations from these two submodules
%and the resulting quantum Lie algebra has the structure described in
%this paper.

The q-conjugation $\sim$ acting on $\uqg$ which we have defined in
definition \ref{qu} does not restrict to $\qlie{g}$. We can however
define a different q-conjugation on $\qlie{g}$.

\begin{definition}
{\underline{q-conjugation on $\qlie{g}$}} is the q-linear map
$\qlie{g}\rightarrow\qlie{g},~a\mapsto a^q$ which extends the
q-conjugation $\sim$ on $\Ch$ by acting as the identity on the
basis elements $X_\a$ and $H_i$.
\end{definition}

The quantum Lie bracket $[a\circ b]$ which we have defined through
the adjoint action in $\uqg$ is clearly not anti-symmetric, i.e.
$[a\circ b]\neq -[b\circ a]$. However we have
\begin{theorem} The quantum Lie bracket is q-anti-symmetric
in the sense that
\[\label{qas}
[a^q\circ b^q]=-[b\circ a]^q,~~~~~\forall a,b\in\qlie{g},~~~g=a_n,c_2.
\]
\end{theorem}
This follows from combining the antiautomorphism $\chi$, described
in equation \Eq{chi},
with the $q$-isomorphism $\t{\th}$: $a^q=-\chi(\t{\th}(a))$.

Our observations in this paper regarding the structure of quantum Lie
algebras have raised many new questions. Among them:
$\bullet$ What is the origin of the q-anti-symmetry \Eq{qas} of the quantum
Lie bracket? This has recently been answered in \ct{Del95b}.
$\bullet$
What are representations of quantum Lie algebras?
$\bullet$ How can the q-symmetric q-bilinear form $\form{\,.\,,\,.\,}$
on root space defined in \Eq{qform} be used to define a q-geometry on
root space? What are q-Weyl "reflections" with respect to such a form?
Can they be used to define quantum root systems axiomatically?
$\bullet$ Is there a connection to quantum affine Toda theory and other
quantum integrable models? These questions are under investigation.

Finally we would like to draw the readers attention to the work of
Sudbery and Lyubashenko 
\cite{Sud2} which has appeared since the completion of this
work. They also give quantum Lie algebras for $sl_2$ and
$sl_3$. For $sl_3$ however they do not impose invariance under the diagram
automorphism.

For further information on quantum Lie algebras visit the quantum
Lie algebra home page on the World Wide Web at 
http://www.mth.kcl.ac.uk/\~delius/q-lie.html.

\vspace{5mm}

\noindent{\bf Acknowledgements}
G.W.D. thanks the Deutsche Forschungsgemeinschaft for a
Habilitationsstipendium. A.H. thanks the EC for a research fellowship.
We have profited greatly from discussions with Mark Gould and Yao-Zhong
Zhang which have led to the work in \ct{Del95}. G.W.D. thanks
Ed Corrigan and Patrick Dorey for discussions and hospitality during
his visit in Durham.

\end{document}